\documentclass[a4paper,12pt]{article}
\usepackage{latexsym}
\usepackage{amssymb,amsmath}
\newcommand{\ux}{\underline{x}}
\newcommand{\uy}{\underline{y}}
\newcommand{\uf}{\underline{f}}
\newcommand{\uh}{\underline{h}}
\newcommand{\ug}{\underline{g}}
\newtheorem{theorem}{Theorem}[section]
\newtheorem{lemma}{Lemma}[section]
\newtheorem{proposition}{Proposition}[section]
\newtheorem{corollary}{Corollary}[section]
\newtheorem{rem}{Remark}[section]
\title{Factorization of the Schr\"{o}dinger Operator and the Riccati Equation in the Clifford Analysis Setting\thanks{This paper was published in \textit{Liber Amicorum Richard Delanghe: een veelzijdig wiskundige}, F. Brackx and H. De Schepper (Eds), Academia Press (Gent), 2005, pp. 69-84.}}
\author{Nele De Schepper\\\normalsize{Department of Mathematical Analysis, Ghent University}\\\normalsize{Galglaan 2, B-9000 Gent, Belgium}\\\normalsize{e-mail: nds@cage.ugent.be}\and
Dixan Pe\~na Pe\~na\\\normalsize{Department of Mathematical Analysis, Ghent University}\\\normalsize{Galglaan 2, B-9000 Gent, Belgium}\\\normalsize{e-mail: dixanpena@gmail.com}}
\begin{document}
\date{}
\maketitle
\section{Introduction}
In one dimension the so-called Riccati equation takes the form
\begin{equation}\label{Riccati}
\frac{dy}{dx} + y^2 = -v \ \ .
\end{equation}
This non-linear equation is also called Miura transformation and gives a relation between the Korteweg-de Vries equation and the modified Korteweg-de Vries equation. It is related to the one-dimensional Schr\"{o}dinger equation
\begin{equation}\label{Schrodinger}
- \frac{d^2u}{dx^2} - v u =0
\end{equation}
with $v$ a function called the potential, in the following way. Given a particular solution for (\ref{Schrodinger}), one can construct the corresponding solution for (\ref{Riccati}) simply by taking the logarithmic derivative $\displaystyle{y=\frac{1}{u} \frac{du}{dx}}$ and vice versa.\\
Moreover, there is a second relation between the one-dimensional Schr\"{o}dinger equation and the Riccati equation, namely the one-dimensional Schr\"{o}dinger operator can be factorized in the form
\begin{equation}\label{factor}
- \frac{d^2}{dx^2} - v = - \left( \frac{d}{dx} + y \right) \left( \frac{d}{dx} - y \right) 
\end{equation}
if and only if $y$ is a solution of (\ref{Riccati}).\\
An analogous factorization of the three-dimensional Schr\"{o}dinger operator is obtained in \cite{Bern2} (see also \cite{Krav1}) using quaternionic analysis, while some representation formulae for the solution of the three-dimensional Schr\"{o}dinger equation are contained in \cite{Bern1}.

The Riccati equation has some peculiar properties among which two theorems of Euler. The first of these states that given a particular solution of the Riccati equation, the general solution can be found in two integrations; while the second one states that given two particular solutions, the general solution can be found in one integration. In \cite{Krav2} a quaternionic generalization of the Riccati equation is established and versions of the above-mentioned Euler theorems corresponding to this generalization are proved.

Closely connected to the Schr\"{o}dinger operator is the Darboux transformation. In general it can be defined as follows. Consider the eigenvalue equation
\begin{equation}\label{eigenv}
Lf=\lambda^2 f 
\end{equation}
with $\lambda$ a complex number. Assume that the operator $L$ admits the factorization $L=AB$ with $B$ a linear operator. Then naturally (\ref{eigenv}) yields
\begin{displaymath}
L_1 f_1 = \lambda^2 f_1 \ \ ,
\end{displaymath}
where $L_1 = BA$ and $f_1=Bf$. Thus applying $B$ to eigenfunctions of $L$ yields eigenfunctions of $L_1$.\\
Although this scheme is general, the Darboux transformation has been studied mostly for the Schr\"{o}dinger operator. In view of the factorization (\ref{factor}), the Schr\"{o}dinger equation
\begin{displaymath}
\left( - \frac{d^2}{dx^2} - v \right) f = \lambda^2 f
\end{displaymath}
has the following Darboux image
\begin{displaymath}
L_1 f_1 = - \left( \frac{d}{dx} - y \right) \left( \frac{d}{dx} + y \right) f_1 = \lambda^2 f_1
\end{displaymath}
where $\displaystyle{f_1 = \left( \frac{d}{dx} -y \right) f}$. The operator $L_1$ has again the form of a Schr\"{o}dinger operator
\begin{displaymath}
L_1= - \frac{d^2}{dx^2} - v_1 \ \ ,
\end{displaymath}
where $\displaystyle{v_1 = \frac{dy}{dx} - y^2}$ .

This paper is devoted to a factorization of the higher dimensional Schr\"{o}- dinger operator in the framework of Clifford analysis, a direct and elegant generalization to higher dimension of the theory of holomorphic functions in the complex plane. This factorization is based on one of the most fundamental features in Clifford analysis, viz. the factorization of the Laplace operator $\Delta_n$ by the so-called Dirac operator $\partial_{\ux}$ : $- \Delta_n = \partial_{\ux}^2$. It is precisely this Dirac operator $\partial_{\ux}$ which underlies the notion of monogenicity of a function, a notion which is the multi-dimensional counterpart to that of holomorphicity in the complex plane.

The outline of the paper is as follows. For the reader who is not familiar with Clifford analysis, we recall some of its basics in Section \ref{CA}. In Section \ref{factorization} we establish a factorization of the higher dimensional Schr\"{o}dinger operator using Clifford analysis. As in the one-dimensional case, this Clifford factorization of the Schr\"{o}dinger operator yields a generalization of the classical Riccati equation to the Clifford analysis setting. Moreover we show that for various reasons this so-called Clifford Riccati equation is a good generalization of the classical one. Some  possibilities for obtaining particular vector-valued solutions of the Clifford Riccati equation are discussed in Section \ref{solRiccati}. Furthermore, in Section \ref{Eulertheorems} generalizations of Euler's theorems for the Clifford Riccati equation are established. Next, a generalized Schr\"{o}dinger operator in Clifford analysis is studied (Section \ref{generalizedSchrodinger}). In a final section a decomposition of the kernel of this generalized Schr\"{o}dinger operator is obtained.
\section{Some basic notions of Clifford analysis}\label{CA}
Clifford analysis (see e.g. \cite{Clif1} and \cite{Clif2}) offers a function theory which is a higher dimensional analogue of the theory of the holomorphic functions of one complex variable. 

The functions considered are defined in the Euclidean space $\mathbb{R}^n$ ($n>1$) and take their values in the Clifford algebra $\mathbb{R}_{0,n}$ or its complexification $\mathbb{C}_n=\mathbb{R}_{0,n} \otimes \mathbb{C}$. If $(e_1,\ldots,e_n) $ is an orthonormal basis of $\mathbb{R}^n$, then a basis for the Clifford algebra $\mathbb{R}_{0,n}$ is given by $ (e_A: A \subset \lbrace 1, \ldots,n \rbrace )$ where $ e_{\emptyset} = 1$ is the identity element. The non-commutative multiplication in the Clifford algebra is governed by the rules:
\begin{displaymath}
e_j e_k + e_k e_j  =  -2 \ \delta_{j,k} \ \ , \quad j,k=1,2,\ldots,n \ \ .
\end{displaymath} 

Conjugation is defined as the anti-involution for which 
\begin{displaymath}
\overline{e_j} = -e_j \ \ , \quad j=1,2,\ldots,n
\end{displaymath}
with the additional rule $\overline{i} = -i$ in the case of $\mathbb{C}_n$ .

For $k=0,1,\ldots,n$ fixed, we call
\begin{displaymath}
\mathbb{R}_{0,n}^{(k)} = \biggl\lbrace a\in\mathbb R_{0,n}:\;a=\sum_{\vert A \vert =k} a_A e_A \ ; \ a_A \in \mathbb{R} \biggr\rbrace
\end{displaymath}
the subspace of $k$-vectors, i.e. the space spanned by the products of $k$ different basis vectors. The 0-vectors and 1-vectors are simply called scalars and vectors respectively, while the $n$-vector $e_N = e_1 e_2 \ldots e_n$ is called the pseudo-scalar.

The Euclidean space $\mathbb{R}^n$ is embedded in the Clifford algebras $\mathbb{R}_{0,n}$ and $\mathbb{C}_n$ by identifying $(x_1,\ldots,x_n)$ with the vector variable $\underline{x}$ given by
\begin{displaymath}
\underline{x} = \sum_{j=1}^n e_j x_j \ \ .
\end{displaymath}
The product of two vectors splits up into a scalar part and a 2-vector or a so-called bivector part:
\begin{displaymath}
\ux \ \uy = \ux \ \bullet \ \uy \ + \ux \wedge \uy \ \ ,
\end{displaymath}
where 
\begin{displaymath}
\ux \ \bullet \ \uy = -\langle \ux, \uy \rangle = - \sum_{j=1}^n x_j y_j
\end{displaymath}
and
\begin{displaymath}
\ux \wedge \uy = \sum_{i=1}^n \sum_{j=i+1}^n e_i e_j (x_i y_j - x_j y_i) \ \ .
\end{displaymath}
Note that the square of a vector variable $\ux$ is scalar-valued and is the norm squared up to a minus sign:
\begin{displaymath}
\ux^2 = - \langle\ux,\ux \rangle = - |\ux|^2 \ \ .
\end{displaymath}

The elliptic vector differential operator of the first order 
\begin{displaymath}
\partial_{\ux} = \sum_{j=1}^n e_j \partial_{x_j}\ \ ,
\end{displaymath}
called Dirac operator, splits the Laplace operator in $\mathbb{R}^n$ :
\begin{displaymath}
\Delta_n = - \partial_{\ux}^2 \ \ .
\end{displaymath}
For a differentiable scalar-valued function $\phi$ and a differentiable Clifford algebra-valued function $f$, we have
\begin{equation}\label{Leibnitz}
\partial_{\ux}( \phi f) = \partial_{\ux}(\phi) \ f + \phi \ \partial_{\ux}(f) \ \ .
\end{equation}
Furthermore, this Leibnitz rule admits the following generalization in the case of a $k$-vector valued function $G_k = \sum_{ \vert A \vert = k} G_{k,A} \ e_A $ (see for example \cite{Dixan})
\begin{equation}\label{Leibnitzspec}
\partial_{\ux}( G_k f) = \partial_{\ux}(G_k) \ f + 2 \sum_{j=1}^n \ \lbrack e_j G_k \rbrack_{k-1} \ \partial_{x_j}(f) + (-1)^k G_k \ \partial_{\ux}(f) 
\end{equation}
where $ \lbrack a \rbrack_k$ denotes the projection of $a \in \mathbb{R}_{0,n}$ on $\mathbb{R}_{0,n}^{(k)}$ .
\section{Clifford factorization of the Schr\"{o}dinger operator and the Clifford Riccati equation}\label{factorization}
In the sequel, $M^{f}$ denotes the operator of multiplication with $f$ from the right, i.e. $M^{f} g = gf$ .\\
The following results generalize to higher dimension the relations mentioned in the introduction between the one-dimensional Schr\"{o}dinger equation and the Riccati equation.
\begin{proposition}
The Schr\"{o}dinger operator may be factorized into
\begin{displaymath}
( - \Delta_n - v I) \phi= \left( \partial_{\ux} + M^f \right) \left( \partial_{\ux} - M^f \right) \phi \ \ ,
\end{displaymath}
where $I$ is the identity operator and $\phi$ a scalar-valued function, if and only if
\begin{equation}\label{CliffRiccati}
\partial_{\ux} f + f^2 = v \ \ .
\end{equation}
\end{proposition}
\textit{Proof.} Applying the Leibnitz rule (\ref{Leibnitz}) yields
\begin{eqnarray*}
\left( \partial_{\ux} + M^f \right) \left( \partial_{\ux} - M^f \right) \phi & = & \left( \partial_{\ux} + M^f \right) (\partial_{\ux} \phi - \phi f)\\
& = & - \Delta_n \phi - \partial_{\ux}(\phi) f - \phi \ \partial_{\ux}(f) + \partial_{\ux}(\phi) f - \phi f^2\\
& = & - \Delta_n \phi - (\partial_{\ux} f + f^2) \phi \ \ .
\end{eqnarray*}
Hence
\begin{displaymath}
( - \Delta_n - v I) \phi= \left( \partial_{\ux} + M^f \right) \left( \partial_{\ux} - M^f \right) \phi
\end{displaymath}
is equivalent with
\begin{displaymath}
\partial_{\ux} f + f^2  = v \ \ .  \ \ \ \ \ \ \square
\end{displaymath}
\begin{proposition}
The scalar-valued function $\phi$ is a solution of the Schr\"{o}dinger equation
\begin{displaymath}
\Delta_n \phi + v \phi = 0
\end{displaymath}
if and only if $\displaystyle{\uf = \frac{\partial_{\ux} \phi}{\phi}}$ is a solution of
\begin{displaymath}
\partial_{\ux} f + f^2 = v \ \ .
\end{displaymath}
\end{proposition}
\textit{Proof.} Assume that $\phi$ is a solution of the Schr\"{o}dinger equation. For $\displaystyle{\uf = \frac{\partial_{\ux} \phi}{\phi}}$ we then have
 \begin{eqnarray*}
 \partial_{\ux} \uf + \uf^2 & = & - \frac{(\partial_{\ux} \phi)^2}{\phi^2} - \frac{\Delta_n \phi}{\phi} + \frac{(\partial_{\ux} \phi)^2}{\phi^2}\\
 & = & v \ \ .
 \end{eqnarray*}
 Conversely, suppose that $\uf$ satisfies $\partial_{\ux} \uf + \uf^2 = v $ and that there exists a function $\phi$ such that $\displaystyle{\uf = \frac{\partial_{\ux} \phi}{\phi}}$. For this function $\phi$ we then obtain
 \begin{displaymath}
 - \frac{(\partial_{\ux} \phi)^2}{\phi^2} - \frac{\Delta_n \phi}{\phi} + \frac{(\partial_{\ux} \phi)^2}{\phi^2} = v
 \end{displaymath}
 or equivalently
 \begin{displaymath}
 \Delta_n \phi + v \phi = 0 \ \ . \ \ \ \ \square
 \end{displaymath}
\\ 
As Propositions 3.1 and 3.2 are the counterparts to the two relations mentioned in the introduction between the Riccati equation (\ref{Riccati}) and the one-dimensional Schr\"{o}dinger equation (\ref{Schrodinger}), equation (\ref{CliffRiccati}) can be considered as a good generalization of (\ref{Riccati}) to the Clifford analysis setting. Moreover, for a one-dimensional solution $f = f(x_k) e_k$ equation (\ref{CliffRiccati}) reduces to the classical one-dimensional Riccati equation
 \begin{displaymath}
 \partial_{x_k}f + f^2 = -v \ \ .
 \end{displaymath}
 Hence we call (\ref{CliffRiccati}) the Clifford Riccati equation.\\
 \begin{rem}
In the special case of a vector-valued solution, the Clifford Riccati equation is equivalent with a scalar elliptic partial differential equation. Indeed, for a vector-valued solution $\uf$ equation (\ref{CliffRiccati}) only consists of the following scalar and bivector part:
 \begin{eqnarray}\label{remark}
 - \langle \partial_{\ux}, \uf \rangle - \langle\uf,\uf \rangle & = & v \ \  ,\\
 \partial_{\ux} \wedge \uf & = & 0 \nonumber \ \ .
 \end{eqnarray}
 The bivector part
 \begin{displaymath}
 \partial_{x_i} f_j - \partial_{x_j} f_i = 0  \quad \mathrm{for} \ i<j 
 \end{displaymath}
 implies that for a simply connected domain $\Omega$ in $\mathbb{R}^n$, there exists a scalar-valued function $\phi$ such that 
 \begin{displaymath}
 \uf = \nabla \phi = \partial_{x_1}(\phi) \ e_1 + \partial_{x_2}(\phi) \ e_2 + \ldots + \partial_{x_n}(\phi) \ e_n \ \ .
 \end{displaymath}
 Inserting the above in the scalar part (\ref{remark}) yields
 \begin{displaymath}
 \Delta_n \phi \ + \langle \nabla \phi, \nabla \phi \rangle = -v \ \ .
 \end{displaymath}
 \end{rem}
 \begin{rem}
 From Proposition 3.2 it is clear that vector-valued solutions of the homogeneous Clifford Riccati equation
 \begin{displaymath}
 \partial_{\ux} \uf + \uf^2 = 0
 \end{displaymath}
 take the form $\displaystyle{\uf = \frac{\partial_{\ux} \phi}{\phi}}$ with $\phi$ a harmonic function, i.e. $\phi \in \mathrm{ker}(\Delta_n)$.
 \end{rem} 
\section{Particular vector-valued solutions of the Clifford Riccati equation}\label{solRiccati}
In this section we discuss some possibilities for obtaining particular vector-valued solutions of the Clifford Riccati equation (\ref{CliffRiccati}).
 
First we consider the special case in which the potential $v$ takes the form
\begin{displaymath}
v(\ux) = v_1(x_1) + v_2(x_2) + \ldots + v_n(x_n) \ \ ,
\end{displaymath}
viz. the variables can be separated. Assuming that
\begin{displaymath} 
\uf = f_1(x_1) e_1 + f_2(x_2) e_2 + \ldots + f_n(x_n) e_n \ \ ,
\end{displaymath}
the Clifford Riccati equation (\ref{CliffRiccati}) reduces to the system of ordinary one-dimensional Riccati equations
\begin{displaymath}
\partial_{x_k} f_k + f_k^2 = - v_k \ \ , \quad k=1,2,\ldots,n \ \ .
\end{displaymath}
Hence in this special case, a particular vector-valued solution of (\ref{CliffRiccati}) can be found if and only if each of the above one-dimensional equations can be solved.

By means of the existence of a large class of vector-valued solutions of the homogeneous Clifford Riccati equation (see Remark 3.2), we are able to reduce the Clifford Riccati equation (\ref{CliffRiccati}) to a scalar differential equation. Indeed, suppose that $\uf$ and $\ug$ are two vector-valued solutions of the homogeneous Clifford Riccati equation, i.e.
\begin{displaymath}
\uf = \frac{\partial_{\ux} \phi_1}{\phi_1} \quad \mathrm{and} \quad \ug = \frac{\partial_{\ux} \phi_2}{\phi_2}
\end{displaymath}
where $\phi_1$ and $\phi_2$ are harmonic. Then the sum $\uf + \ug$ is a solution of (\ref{CliffRiccati}) if and only if we have
\begin{displaymath}
\uf \ \ug + \ug \ \uf = v
\end{displaymath}
or equivalently
\begin{equation}\label{condition}
-2 \left< \frac{\partial_{\ux} \phi_1}{\phi_1} , \frac{\partial_{\ux} \phi_2}{\phi_2} \right> = v \ \ .
\end{equation}
In particular, choosing $\phi_1 = \phi_2 = \phi$, equation (\ref{condition}) becomes
\begin{displaymath}
2 \left(\frac{\partial_{\ux} \phi}{\phi} \right)^2 = v \ \ .
\end{displaymath}
Hence if $\phi \in \mathrm{ker}(\Delta_n)$ is a solution of the above equation, then $ 2 \displaystyle{\frac{\partial_{\ux} \phi}{\phi}}$ is a solution of the Clifford Riccati equation.
 \section{Generalizations of Euler's theorems for the Clifford Riccati equation}\label{Eulertheorems}
 In this section we establish generalizations of Euler's theorems mentioned in the introduction for the Clifford Riccati equation. 
 \begin{proposition}
 Let $\underline h$ be a particular vector-valued solution of the Clifford Riccati equation (\ref{CliffRiccati}). Then the vector-valued function 
 \begin{equation}\label{EulerTheorem1} 
 \underline f=\underline g+\underline h
 \end{equation}
 is also a solution of (\ref{CliffRiccati}) if and only if $\displaystyle{\underline g=\frac{\partial_{\underline x}\phi}{\phi}}$ for some scalar-valued solution $\phi$ of the equation 
 \[\Delta_n\phi+2\langle \partial_{\underline x}\phi,\underline h\rangle=0 \ \ .\] 
 \end{proposition}
 \textit{Proof.} Substituting (\ref{EulerTheorem1}) in (\ref{CliffRiccati}) gives
 \begin{equation} \label{EulerTheorem2}
 \partial_{\underline x}\;\underline g-2\langle \underline g,\underline h\rangle+\underline g^2=0 \ \ .   
 \end{equation}
 As in Remark 3.1, the bivector part of (\ref{EulerTheorem2}) implies that for a simply connected domain $\Omega $ in $\mathbb R^n$, there exists a scalar-valued function $\phi_1$ such that 
 \[\underline g=\partial_{\underline x}\phi_1 \ \ .\] 
 Putting $\phi=\exp(\phi_1)$, we then have that  
 \[\underline g=\frac{\partial_{\underline x}\phi}{\phi}\]
 and equation (\ref{EulerTheorem2}) written in terms of $\phi$ becomes 
 \[- \frac{(\partial_{\ux} \phi)^2}{\phi^2} - \frac{\Delta_n \phi}{\phi} -\frac{2}{\phi}\;\langle \partial_{\underline x}\phi,\underline h\rangle+ \frac{(\partial_{\ux} \phi)^2}{\phi^2} =0\]
 or equivalently
 \[\Delta_n\phi+2\langle \partial_{\underline x}\phi,\underline h\rangle=0 \ \ .\ \ \ \ \square\]
 \begin{proposition}
 Let $\underline g=\partial_{\underline x}\phi_1$, $\underline h=\partial_{\underline x}\phi_2$ be two particular vector-valued solutions of the Clifford Riccati equation (\ref{CliffRiccati}). Then 
 \[\underline f=\frac{\alpha \ \underline g-\underline h}{\alpha-1}\] 
 is also a vector-valued solution of (\ref{CliffRiccati}), where $\alpha=K\exp(\phi_1-\phi_2)$, $K\in\mathbb C$.
 \end{proposition}
 \textit{Proof.} A straightforward computation yields
 \begin{eqnarray*}
 \partial_{\underline x}\underline f+\underline f^2&=&-\frac{\partial_{\underline x}(\alpha)(\alpha\underline g-\underline h)}{(\alpha-1)^2}+\frac{\partial_{\underline x}(\alpha)\underline g+\alpha\partial_{\underline x}\;\underline g-\partial_{\underline x}\;\underline h}{(\alpha-1)}+\frac{(\alpha\underline g-\underline h)^2}{(\alpha-1)^2}\\&=&\frac{-\alpha(\underline g-\underline h)^2+(\alpha\underline g-\underline h)^2+\alpha(\alpha-1)\partial_{\underline x}\;\underline g-(\alpha-1)\partial_{\underline x}\;\underline h}{(\alpha-1)^2}\\&=&\frac{\alpha(\alpha-1)(\partial_{\underline x}\;\underline g+\underline g^2)-(\alpha-1)(\partial_{\underline x}\;\underline h+\underline h^2)}{(\alpha-1)^2}\\&=&v \ \ ,
 \end{eqnarray*}
where we have used the assumption that $\underline g$ and $\underline h$ are particular solutions of the Clifford Riccati equation (\ref{CliffRiccati}).\ \ \ \ $\square$
\begin{rem}
It should be noted that the new solution obtained in Proposition 5.2 is not necessarily the general solution for vector-valued functions. Indeed, let $n>2$ and consider the Clifford Riccati equation with $v=-1$
\[\partial_{\underline x}f+f^2=-1 \ \ .\]  
Two vector-valued solutions of the above equation are $\underline g=e_1=\partial_{\underline x}x_1$ and $\underline h=e_2=\partial_{\underline x}x_2$. Therefore, by Proposition 5.2 we get the class of solutions 
\[\frac{e_1K\exp(x_1-x_2)-e_2}{K\exp(x_1-x_2)-1} \ \ ,\;\;K\in\mathbb C \ \ .\]
However the solutions $e_j$ ($j>2$) are not included in the above expression.
\end{rem}
 \section{A Generalized Schr\"{o}dinger operator}\label{generalizedSchrodinger}
 Let 
 \begin{displaymath}
 \partial_{\ux}^f := \partial_{\ux} + M^f \quad \mathrm{and} \quad \partial_{\ux}^{-f} := \partial_{\ux} - M^f \ \ .
 \end{displaymath}
 In this section we consider the equation
 \begin{equation}\label{genSchrodinger}
 \partial_{\ux}^f \partial_{\ux}^{-f} g = \lambda^2 g
 \end{equation}
 where $g$ is a Clifford algebra-valued function and $\lambda \not= 0$ a complex number. We assume that $f$ is a solution of the Clifford Riccati equation (\ref{CliffRiccati}).\\
From Proposition 3.1 it follows at once that for a scalar-valued function $\phi$ equation (\ref{genSchrodinger}) reduces to the Schr\"{o}dinger equation
 \begin{displaymath}
 (- \Delta_n - v I) \phi = \lambda^2 \phi \ \ .
 \end{displaymath}
 Hence equation (\ref{genSchrodinger}) will be referred to as the generalized Schr\"{o}dinger equation.
 
 Naturally the Darboux scheme described in the introduction works for the generalized Schr\"{o}dinger operator $L = \partial_{\ux}^f \partial_{\ux}^{-f}$.
\begin{proposition}
If $g$ is a solution of (\ref{genSchrodinger}) then $h = \partial_{\ux}^{-f}g$ is a solution of the equation
 \begin{displaymath}
 L_1 h = \lambda^2 h \ \ ,
 \end{displaymath}
 where $L_1 = \partial_{\ux}^{-f} \partial_{\ux}^f$ .\\
 \end{proposition}
 By means of (\ref{Leibnitzspec}) we obtain for a $k$-vector valued function $G_k$ :
 \begin{eqnarray}\label{kvector1}
 & & \partial_{\ux}^{f} \partial_{\ux}^{-f} G_k\nonumber\\
 & = & \left( \partial_{\ux} + M^f \right) \left( \partial_{\ux} G_k - G_k f \right)\nonumber\\
 & = & - \Delta_n G_k - \partial_{\ux}(G_k) f - 2 \sum_{j=1}^{n} \ \lbrack e_j G_k \rbrack_{k-1} \ \partial_{x_j}(f) - (-1)^k G_k \ \partial_{\ux}(f)\nonumber\\
 & &  + \ \partial_{\ux}(G_k) f - G_k f^2\nonumber\\
 & = & - \Delta_n G_k + G_k \left( (-1)^{k+1} \ \partial_{\ux} f - f^2 \right) - 2 \sum_{j=1}^{n}\  \lbrack e_j G_k \rbrack_{k-1} \ \partial_{x_j} f
 \end{eqnarray}
 and similarly
 \begin{equation}\label{kvector2}
 \partial_{\ux}^{-f} \partial_{\ux}^f G_k = - \Delta_n G_k - G_k \left( (-1)^{k+1} \ \partial_{\ux} f + f^2 \right) + 2 \sum_{j=1}^{n} \ \lbrack e_j G_k \rbrack_{k-1} \ \partial_{x_j} f \ \ .
 \end{equation}
 \begin{rem}
 It is easy to check that for a scalar-valued function $\phi$ we have
 \begin{equation}\label{scalar}
 \partial_{\ux}^{-f} \partial_{\ux}^f \phi = - \Delta_n \phi + \phi \ (\partial_{\ux} f - f^2) \ \ .
 \end{equation}
 \end{rem}
Proposition 6.1 and equation (\ref{kvector2}) lead to the following result.
\begin{proposition}
Let $\uf$ be a vector-valued solution of the Clifford Riccati equation and $\phi$ a scalar-valued solution of the Schr\"{o}dinger equation
 \begin{equation}\label{schrodingereq}
 (- \Delta_n - v I) \phi = \lambda^2 \phi \ \ .
 \end{equation}
 Then the vector-valued function
 \begin{displaymath}
 \uh = \sum_{j=1}^n h_j e_j = \partial_{\ux}^{- \uf} \phi
 \end{displaymath}
 is a solution of the equation
 \begin{displaymath}
 \uh (-\Delta_n - v) - 2 \sum_{j=1}^n h_j \ \partial_{x_j} \uf = \lambda^2 \uh \ \ .
 \end{displaymath}
 \end{proposition}
 \textit{Proof.} As (\ref{schrodingereq}) is equivalent with
 \begin{displaymath}
 \partial_{\ux}^{\uf}\  \partial_{\ux}^{-\uf} \ \phi = \lambda^2 \phi \ \ ,
 \end{displaymath}
 Proposition 6.1 implies that the function $\uh = \partial_{\ux}^{-\uf} \ \phi$ is a solution of
 \begin{displaymath}
 \partial_{\ux}^{-\uf} \ \partial_{\ux}^{\uf} \ \uh = \lambda^2 \uh \ \ .
 \end{displaymath}
 In view of (\ref{kvector2}) the above equation is equivalent with
 \begin{displaymath}
 -\Delta_n \uh - \uh (\partial_{\ux} \uf + \uf^2) - 2 \sum_{j=1}^{n} h_j \ \partial_{x_j} \uf  = \lambda^2 \uh
 \end{displaymath}
 which gives the desired result. \ \ \ \ \ \ $\square$\\
 \\
 Similarly we obtain from Proposition 6.1, equations (\ref{kvector1}), (\ref{kvector2}) and (\ref{scalar}) another result.
 \begin{proposition}
Let $\uf$ be a vector-valued solution of the equation 
\[u = \partial_{\ux} \uf - \uf^2.\]  
Furthermore assume that $\ug = \sum_{j=1}^n g_j e_j$ is a solution of 
 \begin{equation}\label{vglg}
 \ug \ (-\Delta_n +u) + 2 \sum_{j=1}^n g_j \ \partial_{x_j} \uf = \lambda^2 \ug \ \ .
 \end{equation}
 Then the scalar-valued function $\phi$ and the bivector valued function $H_2$ constructed in the form
 \begin{displaymath}
 \phi= \lbrack \partial_{\ux}^{-\uf} \ \ug \rbrack_0 \quad , \quad H_2= \lbrack \partial_{\ux}^{-\uf} \ \ug \rbrack_2
 \end{displaymath}
 satisfy the equation
 \begin{displaymath}
 (\phi + H_2) ( - \Delta_n + u) + 2 \sum_{j=1}^n \ \lbrack e_j H_2 \rbrack_1 \ \partial_{x_j} \uf = \lambda^2 (\phi + H_2) \ \ .
 \end{displaymath}
 \end{proposition}
 \textit{Proof.} Equation (\ref{vglg}) can be rewritten as
 \begin{displaymath}
 \partial_{\ux}^{\uf} \ \partial_{\ux}^{-\uf} \ \ug = \lambda^2 \ug \ \ .
 \end{displaymath}
 Applying the Darboux scheme yields that the function $h = \partial_{\ux}^{-\uf} \ \ug$ , which consists only of a scalar and a bivector part, satisfies the equation
 \begin{displaymath}
  \partial_{\ux}^{-\uf} \ \partial_{\ux}^{\uf} \ h = \lambda^2 h \ \ .
 \end{displaymath}
 Denoting by $\phi$ and $H_2$ the scalar, respectively bivector, part of $h$ and using expressions (\ref{kvector2}) and (\ref{scalar}) leads to the desired equation. \ \ \ \ $\square$
\begin{rem} 
Propositions 6.2 and 6.3 are special cases of the following more general result.
 \end{rem}
 \begin{proposition}
Let $\uf$ be a vector-valued solution of the equation 
\[w=(-1)^{k+1} \partial_{\ux} \uf - \uf^2.\]  
Furthermore assume that $G_k$ is a $k$-vector valued solution of
\begin{displaymath}
G_k(- \Delta_n + w) - 2 \sum_{j=1}^{n} \lbrack e_j G_k \rbrack_{k-1} \partial_{x_j} \uf = \lambda^2 G_k \ \ .
\end{displaymath}
Then the $(k-1)$-vector valued function $H_{k-1}$ and the $(k+1)$-vector valued function $H_{k+1}$, constructed in the form
\begin{displaymath}
H_{k-1} = \lbrack \partial_{\ux}^{-\uf} G_k \rbrack_{k-1} \ \ , \quad  H_{k+1} = \lbrack \partial_{\ux}^{-\uf} G_k \rbrack_{k+1}
\end{displaymath}
satisfy the equation
\begin{eqnarray*}
& & (H_{k-1} + H_{k+1}) (-\Delta_n +w) + 2 \sum_{j=1}^{n} \biggl( \lbrack e_j H_{k-1} \rbrack_{k-2} + \lbrack e_j H_{k+1} \rbrack_k \biggr) \partial_{x_j} \uf\\
& = & \lambda^2 (H_{k-1} + H_{k+1}) \ \ .
\end{eqnarray*}
 \end{proposition}
\section{Decomposition of the kernel of the generalized Schr\"{o}dinger operators}\label{deckernel}
We start this section with another representation of the equations (\ref{genSchrodinger}) and
\begin{equation}\label{genSchrodinger2}
\partial_{\ux}^{-f} \partial_{\ux}^f h = \lambda^2 h \ \ , \quad 0 \not= \lambda \in \mathbb{C} \ \ .
\end{equation}
\begin{proposition} 
Let $\displaystyle{\frac{n}{2}}$ be odd and $\uf$ a vector-valued function. Then\\
(i) equation (\ref{genSchrodinger}) can be rewritten in the form $(A^2 - \lambda^2)g = 0,$ where 
\[A=M^{i e_N} \partial_{\ux}^{-\uf} = \partial_{\ux}^{\uf} M^{ie_N};\] 
(ii) equation (\ref{genSchrodinger2}) can be rewritten in the form $(B^2 - \lambda^2)h=0,$ where 
\[B=M^{i e_N} \partial_{\ux}^{\uf} = \partial_{\ux}^{-\uf} M^{i e_N}.\]
\end{proposition}
\textit{Proof.} It is clear that the multiplication operator $M^{ie_N}$ commutes with the Dirac operator $\partial_{\ux}$ . Moreover, the operator $M^{ie_N}$ anti-commutes with the operator $M^{\uf}$. Indeed, as $n$ is even, we have consecutively
\begin{eqnarray*}
M^{\uf} ( M^{ie_N} g) & = & g i e_N \uf\\
& = & (-1)^{n-1} g \uf i e_N\\
& = & - M^{ie_N} (M^{\uf} g ) \ \ .
\end{eqnarray*}
Furthermore, we have
\begin{displaymath}
M^{i e_N} M^{i e_N} = I \ \ ,
\end{displaymath}
since
\begin{eqnarray*}
(ie_N)^2 & = & - e_N^2\\
& = & - (-1)^{n (n+1)/2}\\
& = & 1 \ \ ,
\end{eqnarray*}
where we have used the assumption that $\displaystyle{\frac{n}{2}}$ is odd.\\
In view of the above, we find that
\begin{eqnarray}\label{resu}
M^{ie_N} \partial_{\ux}^{-\uf} M^{ie_N} & = & M^{ie_N} (\partial_{\ux} - M^{\uf} ) M^{ie_N}\nonumber\\
& = & M^{ie_N} M^{ie_N} (\partial_{\ux} + M^{\uf})\nonumber\\
& = & \partial_{\ux}^{\uf} \ \ .
\end{eqnarray}
Statement (i) follows now easily, since
\begin{displaymath}
\partial_{\ux}^{\uf} \partial_{\ux}^{-\uf} = M^{ie_N} \partial_{\ux}^{-\uf} M^{ie_N} \partial_{\ux}^{-\uf} = A^2 \ \ ,
\end{displaymath}
where
\begin{eqnarray*}
A & = & M^{ie_N} \partial_{\ux}^{-\uf}\\
& = & M^{ie_N} (\partial_{\ux} -M^{\uf})\\
& = & (\partial_{\ux} + M^{\uf}) M^{ie_N}\\
& = & \partial_{\ux}^{\uf} M^{ie_N} \ \ .
\end{eqnarray*}
Equation (\ref{resu}) leads in an analogous way to (ii). $ \ \ \ \ \square$\\
\\
Next we prove the following lemma.
\begin{lemma}
Let $\displaystyle{\frac{n}{2}}$ be odd and $\uf$ a vector-valued function. Then one has\\
(i) $\mathrm{ker}(A \pm \lambda) = \mathrm{ker}\left(\partial_{\ux} - M^{\uf\mp \lambda i e_N}\right)$\\
and\\
(ii) $\mathrm{ker}(B \pm \lambda) = \mathrm{ker}\left(\partial_{\ux} + M^{\uf\pm \lambda i e_N}\right)$ .
\end{lemma}
\textit{Proof.} We restrict ourselves to the proof of statement (i), the proof of (ii) being similar.\\
We have consecutively
\begin{eqnarray*}
(A\pm \lambda)g = 0 & \Leftrightarrow & M^{ie_N} (\partial_{\ux} - M^{\uf}) g \pm \lambda g = 0\\
& \Leftrightarrow & (\partial_{\ux} - M^{\uf}) g \pm M^{ie_N} \lambda g = 0\\
& \Leftrightarrow & \partial_{\ux} g - g (\uf \mp \lambda i e_N) = 0\\
& \Leftrightarrow & (\partial_{\ux} - M^{\uf \mp \lambda i e_N}) g = 0 \ \ ,
\end{eqnarray*}
which proves (i) . $ \ \ \  \square$\\
\\
We now arrive at the main result of this section.
\begin{theorem}
Let $\displaystyle{\frac{n}{2}}$ be odd and $\uf$ a vector-valued function. Then the following decompositions hold:\\
(i) $\mathrm{ker}(A^2- \lambda^2) = \mathrm{ker}\left(\partial_{\ux} - M^{\uf+ \lambda i e_N} \right) \oplus \mathrm{ker} \left( \partial_{\ux} - M^{\uf-\lambda i e_N} \right)$.\\
and\\
(ii) $\mathrm{ker}(B^2- \lambda^2) = \mathrm{ker}\left(\partial_{\ux} + M^{\uf - \lambda i e_N} \right) \oplus \mathrm{ker} \left( \partial_{\ux} + M^{\uf+\lambda i e_N} \right)$ .
\end{theorem}
\textit{Proof.} Again we restrict ourselves to the proof of (i). As
\begin{displaymath}
(A^2 - \lambda^2) g = (A-\lambda)(A+\lambda) g = (A+\lambda)(A-\lambda)g \ \ ,
\end{displaymath}
it is clear that
\begin{displaymath}
\mathrm{ker}(A-\lambda) + \mathrm{ker}(A+\lambda) \subset \mathrm{ker}(A^2 - \lambda^2) \ \ .
\end{displaymath}
The converse inclusion 
\begin{displaymath}
\mathrm{ker}(A^2 - \lambda^2) \subset \mathrm{ker}(A-\lambda) + \mathrm{ker}(A+\lambda)
\end{displaymath}
also holds. This can be seen as follows. Take $g \in \mathrm{ker}(A^2 - \lambda^2)$ . By definition we then have
\begin{displaymath}
(A^2 - \lambda^2)g=(A-\lambda) (A+\lambda) g = (A+ \lambda)(A-\lambda)g=0
\end{displaymath}
and consequently
\begin{displaymath}
(A+\lambda)g \in \mathrm{ker}(A-\lambda) \quad \mathrm{and} \quad (A-\lambda)g \in \mathrm{ker}(A + \lambda) \ \ .
\end{displaymath}
Now, decomposing $g$ as:
\begin{displaymath}
g = \frac{1}{2\lambda} (A+\lambda)g - \frac{1}{2 \lambda} (A-\lambda)g \ \ ,
\end{displaymath}
we indeed obtain that $g \in \mathrm{ker}(A-\lambda) + \mathrm{ker}(A+\lambda)$ .\\
Moreover, as $\lambda \not= 0$, we also have
\begin{displaymath}
\mathrm{ker}(A-\lambda) \cap \mathrm{ker}(A+\lambda) = \lbrace 0 \rbrace \ \ ,
\end{displaymath}
hence we have shown that
\begin{displaymath}
\mathrm{ker}(A^2 - \lambda^2) = \mathrm{ker}(A-\lambda) \oplus \mathrm{ker}(A+\lambda) \  \ .
\end{displaymath}
By means of Lemma 7.1 this proves the theorem. $ \ \ \ \ \square$
\\
\begin{corollary}
Let $\displaystyle{\frac{n}{2}}$ be odd and $\uf$ a vector-valued solution of the Clifford Riccati equation. Then any scalar-valued solution $\phi$ of the Schr\"{o}dinger equation
\begin{equation}\label{schrequ}
(-\Delta_n-v) \phi = \lambda^2 \phi
\end{equation}
can be uniquely represented as $\phi = g + h$ , where $g$ and $h$ are Clifford algebra-valued solutions of the equations
\begin{displaymath}
\left( \partial_{\ux} - M^{\uf+\lambda i e_N} \right) g = 0 \quad \mathrm{and} \quad \left( \partial_{\ux} - M^{\uf - \lambda i e_N} \right) h = 0 
\end{displaymath}
respectively.
\end{corollary}
\textit{Proof.} From Proposition 3.1, we know that the Schr\"{o}dinger equation (\ref{schrequ}) is equivalent with 
\begin{equation}\label{equiv}
\partial_{\ux}^{f} \partial_{\ux}^{-f} \phi = \lambda^2 \phi
\end{equation}
where $f$ is a solution of the Clifford Riccati equation. Moreover, taking this solution $f$ vector-valued and $\displaystyle{\frac{n}{2}}$ odd, Proposition 7.1 implies that (\ref{equiv}) can be rewritten in the form
\begin{displaymath}
(A^2-\lambda^2) \phi = 0
\end{displaymath}
hence
\begin{displaymath}
\phi \in \mathrm{ker}(A^2 - \lambda^2) = \mathrm{ker}\left( \partial_{\ux} - M^{\uf+\lambda i e_N} \right) \oplus \mathrm{ker}\left(\partial_{\ux} - M^{\uf-\lambda i e_N} \right) \ \ .
\end{displaymath}
Consequently $\phi$ can be uniquely represented as
\begin{displaymath}
\phi = g+h
\end{displaymath}
with
\begin{displaymath}
g \in \mathrm{ker} \left( \partial_{\ux} - M^{\uf + \lambda  i e_N} \right) \quad \mathrm{and} \quad h \in \mathrm{ker}\left( \partial_{\ux} - M^{\uf - \lambda i e_N} \right) \ \ . \ \ \ \ \square
\end{displaymath} 
\begin{corollary}
Let $\displaystyle{\frac{n}{2}}$ be odd and $\uf$ a vector-valued solution of the equation
\begin{displaymath}
u = \partial_{\ux} \uf - \uf^2 \ \ .
\end{displaymath}
Then any scalar-valued solution $\phi$ of the equation
\begin{equation}\label{coro2}
(- \Delta_n + u) \phi = \lambda^2 \phi
\end{equation}
can be uniquely represented as $\phi=g+h$ , where $g$ and $h$ are Clifford algebra-valued solutions of the equations
\begin{displaymath}
\left(\partial_{\ux} + M^{\uf - \lambda i e_N} \right) g = 0 \quad \mathrm{and} \quad \left( \partial_{\ux} + M^{\uf + \lambda i e_N} \right) h =0
\end{displaymath}
respectively.
\end{corollary}
\textit{Proof.} Using equation (\ref{scalar}) to see that (\ref{coro2}) can be rewritten as
\begin{displaymath}
\partial_{\ux}^{-\uf} \partial_{\ux}^{\uf} \phi = \lambda^2 \phi \ \ ,
\end{displaymath}
the proof runs along the same lines as the proof of the previous corollary. $\square$\ \\
\begin{rem}
If we take
\begin{displaymath}
f = f_1 e_1 + f_2 e_2 + \ldots + f_{n-1} e_{n-1}
\end{displaymath}
instead of a vector-valued $\uf$ and $e_n$ instead of the pseudo-scalar $e_N,$ all results of this section hold for arbitrary dimension, i.e. without the restriction $\displaystyle{\frac{n}{2}}$ odd.
\end{rem}

\end{document}